\documentclass{elsart}
\journal{arXiv.org}
\usepackage{amscd,amsmath,amsfonts,amssymb,xy,graphicx,setspace}
\newcommand{\CC}{\mathbb{C}}
\newcommand{\NN}{\mathbb{N}}
\newcommand{\ZZ}{\mathbb{Z}}

\newtheorem{theorem}{Theorem}
\newtheorem{lemma}{Lemma}

\newtheorem{coro}{Corollary}
\newtheorem{remark}{Remark}

\newcommand{\prf}{\noindent\textbf{Proof.}~}
\newcommand{\eop}{\hfill $\square$\\}
\newcommand{\rep}{\texttt{rep}}
\newcommand{\iss}{\texttt{iss}}
\newcommand{\GL}{\mathsf{GL}}

\newcommand{\Stab}{\mathsf{Stab}}

\newcommand{\db}[1]{\{\!\!\{ #1 \}\!\!\}}

\newcommand{\Der}{\mathrm{Der}}
\renewcommand{\sl}{\mathfrak{sl}}
\renewcommand{\sp}{\mathfrak{sp}}
\newcommand{\h}{\mathfrak{h}}
\newcommand{\nl}{\mathfrak{n}}

\newcommand{\cO}{\mathcal{O}}
\newcommand{\cN}{\mathcal{N}}
\newcommand{\cT}{\mathcal{T}}

\xyoption{all}
\xyoption{poly}
\allowdisplaybreaks

\begin{document}
\begin{frontmatter}
\title{On the Structure of the Necklace Lie Algebra}
\author{Jacques Alev}
\ead{jacques.alev@univ-reims.fr}
\address{
Laboratoire de Math\'ematiques\\
U.F.R. Sciences Exactes et Naturelles\\
Universit\'e de Reims\\
Moulin de la Housse - BP 1039\\ 
F-51687 REIMS cedex 2 (France)\\
}
\author{Geert Van de Weyer\thanksref{marie}}
\ead{geert.vandeweyer@ua.ac.be}
\address{
Laboratoire de Math\'ematiques\\
U.F.R. Sciences Exactes et Naturelles\\
Universit\'e de Reims\\
Moulin de la Housse - BP 1039\\ 
F-51687 REIMS cedex 2 (France)\\
}
\thanks[marie]{
The second author was supported by a grant in the Marie Curie Research Training Network MRTN-CT 2003-505078.}

\begin{abstract}
The necklace Lie algebra for a quiver was introduced simultaneously by Bocklandt and Le Bruyn in \cite{RafLieven} and Ginzburg in \cite{Ginz} to study the noncommutative symplectic geometry of preprojective algebras. Later, Van den Bergh in \cite{DPA} linked the necklace Lie algebra to double Poisson structures and used it to define the structure of Poisson varieties on the spaces of isomorphism classes of semi simple representations of an algebra $A$.

In this note, we initiate the systematic study of the Lie algebra structure of the necklace Lie algebra $\nl$ of a free algebra in $2d$ variables. We begin by giving a description of $\nl$ as an $\sp(2d)$-module. Specializing to $d = 1$, we decompose $\nl$ into a direct sum of highest weight modules for $\sl_2$, the coefficients of which are given by a closed formula. Next, we observe that $\nl$ has a nontrivial center, which we link through the center $C$ of the trace ring of couples of generic $2\times 2$ matrices to the Poisson  center of $S(\sl_2)$. The Lie algebra structure of $\nl$ induces a Poisson structure on $C$, the symplectic leaves of which we are able to describe as coadjoint orbits for the Lie group of the semidirect product $\sl_2\rtimes\h$ of $\sl_2$ with the Heisenberg Lie algebra $\h$. Finally, we provide a link between double Poisson algebras on one hand and Poisson orders, introduced by Brown and Gordon in \cite{PoissonOrder} on the other hand, showing that all trace rings of a double Poisson algebra are Poisson orders over their center.
\end{abstract}

\end{frontmatter}

\section{Introduction}
In the last decade there has been an increasing interest in the study of noncommutative symplectic geometry, motivated amongst others by the study of deformed preprojective algebras in general and Calogero-Moser phase space in particular. An important concept introduced in this subject is the notion of the necklace Lie algebra, which was used by Ginzburg in \cite{Ginz} to show that the space of isomorphism classes of semi simple representations of a preprojective algebra is very similar to a coadjoint orbit. Later, Van den Bergh in \cite{DPA} placed the necklace Lie algebra in a more general context of \emph{noncommutative Poisson geometry}, and it is this more general point of view which motivates the study of the necklace Lie algebra in this note.

When one wants to define a meaningful framework of noncommutative algebraic geometry, there are several approaches possible. One approach, advocated by Kontsevich in e.g.~\cite{KontAppr} and Le Bruyn in \cite{LievenBook}, is by \emph{approximation} and is motivated by the observation that for any algebra $A$, we have an infinite number of classical algebraic varieties associated to it, namely its representation spaces. In particular, if $A$ is an affine associative algebra, that is, $A$ is given by a finite number of generators subject to the relations in an ideal $I$,
$$A = \frac{\CC\langle x_1,\dots,x_d\rangle}{I},$$
we have the spaces
$$\rep_n(A) = \{R\in M_n(\CC)^{\oplus d}\mid\forall f\in I: f(R) = 0\}.$$
Each of these spaces has a natural action of the general linear reductive group $\GL_n(\CC)$ on them by conjugation. That is, for $g\in\GL_n(\CC)$ and $R = (R_1,\dots,R_d)\in\rep_n(A)$ we have
$$g.R := (gR_1g^{-1},\dots,gR_dg^{-1}).$$
The quotient space with respect to this action is denoted by $\iss_n(A)$ and classifies isomorphism classes of semi simple representations of $A$. The main idea of the \emph{approximation strategy} now states that a suitable notion in noncommutative geometry should be visible on \emph{all} representation spaces. The archetypical example of such a noncommutative geometric notion is a \emph{formally smooth algebra} \cite{CuntzQuillen}. Recall that an affine associative algebra is called formally smooth if it satisfies the following lifting property. For all affine associative algebras $B$, nilpotent ideals $I$ in $B$ and maps $\varphi:A\rightarrow B/I$ we have
$$\xymatrix{A\ar[dr]^{\varphi}\ar@{-->}[r]^{\exists\tilde{\varphi}} & B \ar@{->>}[d] \\ & B/I.
}$$
For a formally smooth algebra $A$, it is known that all its finite dimensional representation spaces $\rep_n(A)$ are \emph{smooth varieties}.

As all representation spaces $\rep_n(A)$ have natural actions of the corresponding $\GL_n(\CC)$ on them, it is also meaningful in this context to determine what one can say about the structure of all quotient spaces $\iss_n(A)$. Whereas the geometry of all $\rep_n(A)$ is governed by additional structure on $A$, the geometry of all $\iss_n(A)$ is governed by additional structure on the \emph{necklaces of $A$}. That is, the geometry of $\iss_n(A)$ is governed by additional structure on the vector space quotient
$$\nl_A := \frac{A}{[A,A]}.$$
This follows immediately from Procesi's Theorem on $\GL_n$-invariants, which states that the coordinate ring $\CC[\iss_n(A)] = \CC[\rep_n(A)]^{\GL_n(\CC)}$ of $\iss_n(A)$ is generated by traces of words in generic matrices $X_1,\dots,X_d$. That is, we know the coordinate ring $\CC[\rep_n(A)]$ is generated by elements $x_{i,kl}$ which simply project $R\in\rep_n(A)$ to the $(k,l)$-th entry in $R_i$. Now consider the \emph{generic matrices} $X_i = (x_{i,kl})$, then $\CC[\iss_n(A)]$ is generated by elements of the form $tr(X_{i_1}\dots X_{i_m})$ with $m\in\NN$, where $tr$ is the \emph{trace map}. Now note that this is actually the image of the element $x_{i_1}\dots x_{i_m}\in A$ under the map which replaces each generator of $A$ by its corresponding generic matrix and then takes the trace of the matrix thus obtained. Finally, observe that the trace map is a linear map which is invariant under cyclic permutation, so the map considered is actually a map from the vector space quotient $\nl_A$ to $\CC[\iss_n(A)]$, from which we deduce that it is indeed $\nl_A$ which should be studied if one wants to make statements about all quotient spaces at the same time.

\begin{figure}\label{necklace}
\[
\begin{xy}/r3pc/:
{\xypolygon10{~*{\bullet}~>>{}}},
"1" *+{\txt{\tiny $u$}},"8"="a","10"="a1","9"="a2","0"="c1",
"c1"+(0,-1.9),
{\xypolygon10{~*{\bullet}~>>{}}},
"5" *+{\txt{\tiny $v$}},"4"="b","6"="b1","5"="b2","0"="c2",
"a";"a2" **@{/};
"a" +(-2,0) *+{\txt{$\underset{a \in Q_1}{\sum}$}};
"c1" *+{\txt{$\overset{\curvearrowleft}{w_1}$}};
"c2" *+{\txt{$\overset{\curvearrowleft}{w_2}$}};
\POS"a" \ar@/^2ex/^{\txt{\tiny{$a$}}} "a2"
\POS"a2" \ar@/^2ex/^{\txt{\tiny{$a^*$}}} "a"
\end{xy}~\qquad
   \begin{xy}/r3pc/:
{\xypolygon10{~*{\bullet}~>>{}}},
"1" *+{\txt{\tiny $u$}},"8"="a","10"="a1","9"="a2","0"="c1",
"c1"+(0,-1.9),
{\xypolygon10{~*{\bullet}~>>{}}},
"5" *+{\txt{\tiny $v$}},"4"="b","6"="b1","5"="b2","0"="c2",
"a";"a2" **@{/};
"a"+(-2,0) *+{-};
"c1" *+{\txt{$\overset{\curvearrowleft}{w_2}$}};
"c2" *+{\txt{$\overset{\curvearrowleft}{w_1}$}};
\POS"a" \ar@/^2ex/^{\txt{\tiny{$a$}}} "a2"
\POS"a2" \ar@/^2ex/^{\txt{\tiny{$a^*$}}} "a"
\end{xy}
\]
\caption{Lie bracket $[ w_1,w_2 ]$ in $\nl_{\CC\overline{Q}}$. For each occurrence of an arrow $a\in Q$ in $w_1$ we look for each occurrence of $a^*$ in $w_2$, we remove $a$ from $w_1$ and $a^*$ from $w_2$ and connect the opened necklaces. We add up all necklaces thus obtained and then reverse the roles of $w_1$ and $w_2$, subtracting this from the first expression.}
\label{bracket}
\end{figure}
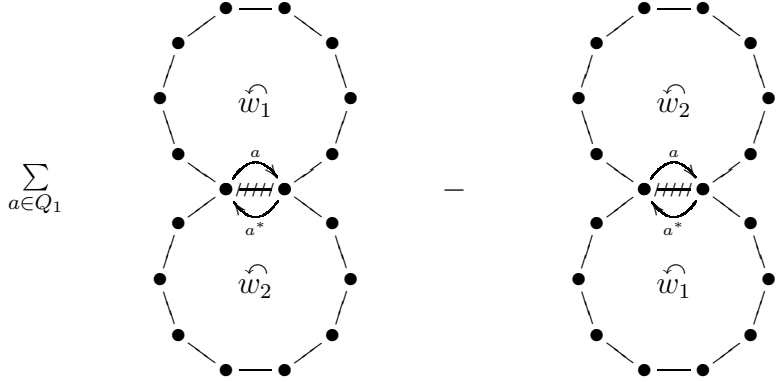

Suppose in particular that $A = \CC \overline{Q}$ is the path algebra of the \emph{double of a quiver}. That is, we have an oriented graph $Q$ and we turn this into a symmetric oriented graph $\overline{Q}$ by adding for each arrow $a$ in $Q$ an arrow $a^*$ in the opposite direction. Then it was observed by Kontsevich in \cite{KonSym} that the necklaces in $\nl_{\CC\overline{Q}}$ have a natural structure of a Lie algebra, which is called the necklace Lie algebra. The definition of this \emph{Kontsevich Lie bracket} is recalled in Figure \ref{necklace}.
In \cite{DPA}, Van den Bergh now observes that if $A$ is a \emph{double Poisson algebra}, then its necklaces have the structure of a Lie algebra. Moreover, the corresponding Lie-Poisson structure on $S(\nl_A)$ induces Poisson structures on all quotient spaces $\iss_n(A)$.

Motivated by this last observation, we initiate the systematic study of the most canonical example of a necklace Lie algebra, namely the necklace Lie algebra of the free algebra $A = \CC\langle x_1,\dots,x_d,x_1^*,\dots,x_d^*\rangle$ in $2d$ variables, equipped with the Kontsevich Lie bracket. Alternatively, this can be seen as the necklace Lie algebra of the free algebra $A$, equipped with the \emph{canonical symplectic double bracket} defined on the generators of $A$ as
$$\db{x_i,x_j} = 0, \db{x_i^*,x_j^*} = 0\quad\mathrm{for~all~} i,j,$$
$$\db{x_k,x^*_\ell} = 0 \quad\mathrm{for~all}\quad k\neq l \quad\mathrm{and}\quad \db{x_i,x_i^*} = 1\otimes 1.$$
This necklace Lie algebra, which from now on we will denote by $\nl$, then induces a Poisson bracket on the center $\CC[\iss_n(A)]$ of the trace ring of $2d$-tuples of $n\times n$ generic matrices. Observe that this necklace Lie algebra can be considered as a generalization of the symplectic Lie algebra, and this symplectic viewpoint was already explored by Crawley-Boevey, Etingof and Ginzburg in \cite{CBEG}, where in particular the center of this Lie algebra is determined.

The necklace Lie algebra carries a natural grading induced from the natural grading of $A$, and we begin our study by observing that homogeneous parts in lowest degree have a nice description in terms of classical Lie algebras. In particular, we have that $\nl_{\leq 1} \cong \h(2d+1)$,
$\nl_2 \cong \sp(2d)$ and $\nl_{\leq 2}\cong \h(2d+1)\rtimes\sp(2d)$. The fact that $\nl_2 \cong \sp(2d)$ in combination with the fact that the necklace Lie bracket is graded of degree $-2$ then allows us to consider $\nl$ as an $\sp(2d)$-module, which is seen to be canonically identified with a quotient of the tensor representation of the canonical representation of $\sp(2d)$ on 
$$V = \CC x_1\oplus \dots\oplus \CC x_d \oplus \CC x_1^*\oplus \dots\oplus \CC x_d^*.$$
Specializing to $d = 1$, we are able to give an explicit Clebsch-Gordan decomposition of $\nl$ as an $\sl_2$-module, where the homogeneous part of degree $n$, $\nl_n$ decomposes in highest weight modules $\nl_{n,m}$ of weight $n-2m$ which occur with multiplicity
\begin{align*}
{n\choose m} - {n \choose m-1} &
- \sum_{\ell |n | \ell m} \frac{\ell-1}{\ell}\sum_{k |\gcd(\ell,\frac{\ell}{n}m)} \mu(k){\frac{\ell}{k} \choose \frac{\ell m}{kn}} \\
&+  \sum_{\ell |n | \ell (m-1)} \frac{\ell-1}{\ell}\sum_{k |\gcd(\ell,\frac{\ell}{n}(m-1))} \mu(k){\frac{\ell}{k} \choose \frac{\ell (m-1)}{kn}},
\end{align*}
where $\mu$ is the M\"obius function.
Next, we observe that for $d = 1$ the necklace Lie algebra has an infinite dimensional center, which means in particular that as Poisson varieties, all $\iss_n(\CC\langle x,x^*\rangle)$ have nontrivial Poisson center, implying for instance that they all have infinite dimensional zero-th Poisson cohomology.

We conclude by observing that for a general double Poisson algebra $A$, for each $n\in\NN$ the trace ring of generic matrices of $A$ is a Poisson order, linking the notion of a double Poisson algebra to another notion of noncommutative Poisson geometry introduced in the study of symplectic reflection algebras.

We organized the paper as follows. In Section \ref{prelims}, we gather all preliminary material needed for the remainder of the paper. In Section \ref{modstructure} we show the homogeneous part of degree two of the necklace Lie algebra is isomorphic to $\sp(2d)$ and then use this description for $d = 1$ to decompose $\nl$ as $\sl(2)$-module. In Section \ref{center} we then turn our attention towards the center $\mathfrak{c}$ of the necklace Lie algebra, showing it is infinite dimensional and linking this center with the unique generating Casimir of the symmetric algebra of $\sl(2)$.  We then turn our attention to the lowest dimensional approximations of the noncommutative geometry of $\CC\langle x,x^*\rangle$, and describe the symplectic leaves of $\iss_2(\CC\langle x,x^*\rangle)$ as coadjoint orbits in Section \ref{leaves}. We conclude the paper with a link between double Poisson algebras and Poisson orders in Section \ref{order}.

Throughout the paper, all associative algebras will be unital over $\CC$. Unadorned tensor products will be tensor products over the base field $\CC$.

\section{Preliminaries}\label{prelims}
In this section, we gather the necessary preliminary material needed for the remainder of the paper. We refer to \cite{DPA} for the details of the results summarized here.

We will consider two actions of an affine associative $\CC$-algebra $A$ on $A\otimes A$. The \emph{outer action} of $A$ on $A\otimes A$ is defined as
$$a.(b'\otimes b'').c := (ab')\otimes (b''c),$$
whereas the \emph{inner action} is defined as
$$a\circ(b'\otimes b'')\circ c := (b'c)\otimes(ab'')$$
for all $a,b',b''\in A$. 

There is a canonical action of the cyclic permutation $\sigma = (1\dots n)$ on elements of $A^{\otimes n}$, defined as
$$\sigma.a_1\otimes\dots\otimes a_n := a_n\otimes a_1\otimes\dots\otimes a_{n-1},$$
which extends to a canonical action of the cyclic group $\ZZ_n$ generated by $\sigma$.
We will use the following shorthand notation for the action of $\sigma$ and $\sigma^{-1}$:
$$\overrightarrow{w} := \sigma.w \quad\mathrm{and}\quad \overleftarrow{w} := \sigma^{-1}.w.$$
Now recall
\begin{defn}
Let $A$ be an affine associative $\CC$-algebra, then a linear map
$$\db{-,-} : A\otimes A\rightarrow A\otimes A$$
is called a \emph{double Poisson bracket} if and only if it satisfies
\begin{enumerate}
\item $\forall a,b\in A : \db{a,b} = - \db{b,a}^o$;
\item $\forall a,b,c\in A : \db{a,bc} = b.\db{a,c} + \db{a,b}.c$;
\item $\forall a,b,c\in A:$
$$\db{a,\db{b,c}'}\otimes \db{b,c}'' + \overrightarrow{\db{b,\db{c,a}'}\otimes \db{c,a}''} + \overleftarrow{\db{c,\db{a,b}'}\otimes \db{a,b}''} = 0.$$
\end{enumerate}
The last identity is known as the \emph{double Jacobi identity}.

An algebra equipped with a double Poisson bracket is called a \emph{double Poisson algebra}.
\end{defn}
In \cite{DPA} it was observed that if $A$ is a double Poisson algebra with double Poisson bracket $\db{-,-}$, this double bracket induces something very close to a Poisson bracket on $A$. Recall that
\begin{defn}
A \emph{Loday algebra} is a vector space $L$ together with a bilinear map
$$\{-,-\}:L\times L \rightarrow L$$
which satisfies :
$$\forall a,b,c\in L : \{a,\{b,c\}\} = \{\{a,b\},c\} + \{b,\{a,c\}\}.$$
A \emph{Loday-Poisson algebra} is an affine associative algebra $L$ together with a bilinear map $\{-,-\}:L\times L \rightarrow L$ which is a derivation in its second argument and which makes $L$ into a Loday algebra.
\end{defn}
We then have
\begin{theorem}
If $A$ is a double Poisson algebra with double bracket $\db{-,-}$ and multiplication map $\mu_A$, then  the associated bracket
$$\{-,-\}_L : A\otimes A \rightarrow A : a\otimes b \mapsto \mu_A(\db{a,b})$$
\begin{enumerate}
\item defines a Loday-Poisson algebra structure on $A$ which satisfies
$$\forall a,b,c\in A: \{[a,b],c\} = 0;$$
\item induces a Lie algebra structure, denoted by $\{-,-\}$ on the vector space quotient
$$\nl_A := \frac{A}{[A,A]}$$
of the algebra $A$ by the vector space spanned by all commutators in $A$.
\end{enumerate}
\end{theorem}
If $A$ is a double Poisson algebra, we will call $\nl_A$ with the induced Lie algebra structure the \emph{necklace Lie algebra of $A$}. 

In general, $\nl_A$ is an infinite dimensional Lie algebra which is generated as a vector space by the equivalence classes of words in the generators of $A$ under cyclic permutation. These basis elements will be called \emph{necklaces}. Observe there is a natural grading on the necklace Lie algebra by the length of the necklaces. This necklace Lie algebra structure now induces Poisson structures on all spaces $\iss_n(A)$:
\begin{theorem}\label{inducedstructure}
Let $A$ be a double Poisson algebra with corresponding necklace Lie algebra $\nl_A$ and consider the induced Lie-Poisson structure on the symmetric algebra $S(\nl_A)$, then the kernel of the canonical algebra morphism
$$S(\nl_A)\twoheadrightarrow \CC[\iss_n(A)],$$
induced by the linear map $\nl_A \rightarrow\CC[\iss_n(A)] : w\mapsto tr(w)$, is a Poisson ideal and hence $\CC[\iss_n(A)]$ is a Poisson algebra.
\end{theorem}

Now consider the free algebra $A = \CC\langle x_1,\dots, x_d,x_1^*,\dots, x_d^*\rangle$ then this algebra has a canonical double Poisson bracket on it defined as
\begin{eqnarray*}
\forall 1\leq i,j\leq d &:& \db{x_i,x_j} = \db{x_i^*, x_j^*} = 0\\
\forall 1\leq i\neq j\leq d &:& \db{x_i,x_j^*} = 0\\
\forall 1\leq i \leq n &:& \db{x_i,x_i^*} = 1\otimes 1
\end{eqnarray*}
on the generators, and extended to a double Poisson bracket using the fact that the bracket has to be twisted antisymmetric and a derivation with respect to the outer bimodule structure. Observe that this bracket is a graded bracket of degree $-2$ for the canonical grading of $A$. In particular, this bracket respects the natural grading of the necklace Lie algebra.

For the greatest part of the remainder of this paper, we will study this particular necklace Lie algebra, which we will simply denote by $\nl$. We will consider it as a Lie algebra with an underlying graded vector space and a bracket of degree $-2$.

\section{The $\sl_2$-Module Structure of $\nl$.}\label{modstructure}
We begin by observing there is a nice description of the Hilbert series of the necklace Lie algebra.
\begin{prop}
Let $\nl$ be the necklace Lie algebra in $2d$ variables and let $\nl_k$ be its homogeneous part of degree $k$, then
$$\dim_\CC \nl_k = \frac{1}{k}\sum_{i=1}^{k} (2d)^{\gcd(k,i)}.$$
\end{prop}
\prf
This is a well-known expression for the number of different necklaces with $2d$ beads, which is a version of \emph{Moreau's necklace counting function} \cite{Moreau}.
\eop

Next, observe that
\begin{theorem}\label{LieIso}
Let $\nl$ be the necklace Lie algebra in $2d$ variables, and let $A_{ab} = \CC[x_1,\dots,x_d,x_1^*,\dots,x_d^*]$ be the polynomial algebra in $2d$ variables equipped with its canonical symplectic Poisson bracket, then 
\begin{enumerate}
\item all homogeneous subspaces of $\nl$ are $\nl_2$-modules;
\item we have a graded morphism of Lie algebras
$$\varphi:\nl\rightarrow A_{ab} : w\mapsto w \mod ([A,A]);$$
\item in particular, there are Lie algebra isomorphisms
\begin{enumerate}
\item $\nl_{\leq 1} \cong \h(2d+1),$
\item $\nl_2 \cong \sp(2d),$
\item $\nl_{\leq 2}\cong \h(2d+1)\rtimes\sp(2d)$.
\end{enumerate}
\end{enumerate}
\end{theorem}
\prf
The fact that $\nl_2$, $\nl_{\leq 1}$ and $\nl_{\leq 2}$ are Lie subalgebras of $\nl$ as well as the fact that all homogeneous subspaces $\nl_i$, $i\in\NN$, are $\nl_2$ modules follows immediately from the fact that the bracket on the necklace Lie algebra is of degree $-2$.

Now observe that by Theorem \ref{inducedstructure}, we have a Poisson algebra morphism
$$S(\nl)\twoheadrightarrow \CC[\iss_1(A)] = A_{ab} =  \CC[x_1,\dots,x_d,x_1^*,\dots,x_d^*],$$
which maps a generator $x_i$ to $tr(x_i)$. Because we are looking at one-dimensional representations, $tr(x_i)$ is canonically identified with $x_i$, so the Poisson bracket on $\CC[x_1,\dots,x_d,x_1^*,\dots,x_d^*]$, defined by the necklace Lie algebra becomes
\begin{eqnarray*}
\forall 1\leq i,j\leq d &:& \{x_i,x_j\}_{ab} = \{x_i^*, x_j^*\}_{ab} = 0\\
\forall 1\leq i\neq j\leq d &:& \{x_i,x_j^*\}_{ab} = 0\\
\forall 1\leq i \leq n &:& \{x_i,x_i^*\}_{ab} = \mu(1\otimes 1) = 1,
\end{eqnarray*}
which exactly corresponds to the canonical symplectic Poisson bracket on $A_{ab}$.
But then this morphism induces in particular a graded morphism of Lie algebras
$$\varphi: \nl\rightarrow A_{ab} : w\mapsto w \mod ([A,A]).$$
Now observe that $\varphi$ is an isomorphism when restricted to $\nl_2$, as both $\nl_2$ and $A_{ab,2}$ are generated by all necklaces of the form $xy$ with
$$x,y\in\{x_1,\dots,x_d,x_1^*,\dots,x_d^*\}.$$
Finally, note that the symplectic Lie algebra $\sp(2d)$ is the restriction of $\{-,-\}_{ab}$ to the degree two subspace of $A_{ab}$, whence the desired isomorphism.

Similarly, $\varphi$ is an isomorphism when restricted to $\nl_{\leq 1}$ and the Heisenberg Lie algebra $\h(2d+1)$ is easily seen to be the restriction of $\{-,-\}_{ab}$ to $A_{ab,\leq 1}$. Finally, as $\nl_{\leq 2} = \nl_{\leq 1}\oplus \nl_2$ is isomorphic as a graded Lie algebra to $A_{ab,\leq 2} = A_{ab,\leq 1}\oplus A_{ab,2}$, and as the bracket on $A_{ab,\leq 2}$ corresponds to the semidirect product, we have the last isomorphism.
\eop
An immediate corollary from this theorem is
\begin{coro}\label{notsimplealgebra}
Denote by $\mathfrak{c}$ the center of $\nl$.
We have the following strict inclusions of infinite dimensional Lie algebras
$$0 \subsetneq \mathfrak{c}\subsetneq \ker(\varphi) \subsetneq \nl_{\geq 4}$$
\end{coro}
\prf
From Theorem \ref{LieIso} we know we have a graded morphism of Lie algebras from $\nl$ to $A_{ab}$. It now simply suffices to observe that the dimensions of the homogenous parts of degree $n$ of the necklace Lie algebra are strictly bigger than the dimensions of the corresponding homogeneous parts of $A_{ab}$ as soon as $n\geq 4$.

Moreover, the fact that $\mathfrak{c}$ is included in $\ker(\varphi)$ follows from
\cite[Cor. 8.6.2]{CBEG}, which states that $\mathfrak{c}$ is generated by all powers $c^n$, of the element $$c := \sum_{i=1}^d [x_i,x_i^*].$$
\eop

Next, let us take a closer look at the action of $\nl_2$ on the homogeneous subspaces $\nl_i$. Consider the vector space
$$V = \CC x_1\oplus\dots\oplus\CC x_d\oplus\CC x_1^*\oplus\dots\oplus\CC x_d^*,$$
and consider the natural action of $\sp(2d)$ on $V$ induced from the canonical Poisson bracket on $\CC[V]$. Recall that the canonical extension to an action on the tensor power $V^{\otimes n}$ is given by
$$g.v_1\otimes\dots\otimes v_n := \sum_{i=1}^n v_1\otimes\dots\otimes v_{i-1}\otimes g.v_i\otimes v_{i+1}\otimes\dots\otimes v_n.$$
\begin{lemma}
Let
$$C_n := \langle v_1\otimes\dots\otimes v_n - v_2\otimes\dots\otimes v_n\otimes v_1 \mid v_1,\dots, v_n\in V\rangle,$$
then $C_n$ is an $\sp(2d)$-submodule of $V^{\otimes n}$.
\end{lemma}
\prf
Indeed, we have for $g\in\sp(2d)$ that
\begin{align*}
&g.(v_1\otimes\dots\otimes v_n - v_2\otimes\dots\otimes v_n\otimes v_1) \\
=&\sum_{i=1}^n v_1\otimes\dots\otimes v_{i-1}\otimes g.v_i\otimes v_{i+1}\otimes\dots\otimes v_n \\
&- \sum_{j=2}^{n-1} v_2\otimes\dots\otimes v_{i-1}\otimes g.v_i\otimes v_{i+1}\otimes\dots\otimes v_n\otimes v_1 - v_2\otimes\dots\otimes v_n\otimes g.v_1 \\
=& \sum_{i=2}^{n-1} (v_1\otimes\dots\otimes v_{i-1}\otimes g.v_i\otimes v_{i+1}\otimes\dots\otimes v_n \\
& - v_2\otimes\dots\otimes v_{i-1}\otimes g.v_i\otimes v_{i+1}\otimes\dots\otimes v_n\otimes v_1) \\
  & +g.v_1\otimes\dots\otimes v_n- v_2\otimes\dots\otimes v_n\otimes g.v_1\\
  & \in C(n).
\end{align*}
\eop

So we have a canonical action of $\sp(2d)$ on $V^c_n := V^{\otimes n}/C(n) = \nl_n$. In particular, we have
\begin{theorem}
Let $\nl$ be the necklace Lie algebra in $2d$ variables and consider $\nl_n$ as a $\sp(2d)$-module through the canonical identification $\nl_2 =  \varphi^{-1}(\sp(2d))$, then the canonical identification 
$$V_n^c \tilde{\rightarrow} \nl_n : v_1\otimes\dots\otimes v_n\mapsto v_1\dots v_n$$
is an isomorphism of $\sp(2d)$-modules.
\end{theorem}
\prf
Recall that the Lie algebra structure on $\nl$ is induced by the Loday-Poisson bracket on $A$ through
$$\{w_1,w_2\} = \{\tilde{w}_1,\tilde{w}_2\}_L \mod [A,A]$$
for necklaces $w_1$ and $w_2$. Now $\{-,-\}_L$ is a derivation in its second argument, so we have for $w\in \nl_2$ and $v = v_1\dots v_n\in \nl_n$ that
$$w.v = \{w,v_1\dots v_n\} = \sum_{i = 1}^n v_1\dots v_{i-1}\{w,v_i\}_Lv_{i+1}\dots v_n \mod [A,A],$$
which exactly corresponds to the image of
$$w.(v_1\otimes\dots\otimes v_n) = \sum_{i=1}^n v_1\otimes\dots\otimes v_{i-1}\otimes\{w,v_i\}\otimes v_{i+1}\otimes\dots\otimes v_n.$$
\eop

From now on, assume $\nl$ is the necklace Lie algebra in $2$ variables, so $\nl_2$ becomes isomorphic to $\sl(2)$, where the isomorphism identifies
$$E = \frac{(x^*)^2}{2},\quad F = -\frac{x^2}{2}\quad\mathrm{and}\quad H = xx^*,$$
with
$$[H,E] = 2E,\quad [H,F] = -2F \quad\mathrm{and}\quad [E,F] = H.$$
The graded Lie algebra morphism from Theorem \ref{LieIso} becomes a morphism from $\nl_2$ to $\CC[x,x^*]$. In this commutative setting, all homogeneous spaces $\CC[x,x^*]_n$ become simple $\sl_2$-modules. However, their noncommutative counterparts in the necklace Lie algebra will in general be no longer simple. We have
\begin{prop}
Let $\nl$ be the necklace Lie algebra in two variables. If $n\geq 4$ then $\nl_n$ is not a simple $\sl_2$ representation.
\end{prop}
\prf
This follows immediately from the fact that the morphism $\varphi_n$ in particular is a morphism of $\sl_2$-modules which is an isomorphism if and only if $n\leq 3$, as can be seen easily by looking at the dimensions of the corresponding homogeneous subspaces.
\eop

If we want to describe $\nl$ as a $\sl_2$-module, we first have to study the $\sl_2$-module structure of $C_n$. Of particular use in this study is the action of the generator $H = xx^*$ on words in $V^{\otimes n}$. We have
\begin{lemma}\label{eigenvector}
Let $w = w_1\otimes\dots\otimes w_n\in V^{\otimes n}$, denote the total degree of $w$ in $x$ as $\deg_x(w)$ and the total degree of $w$ in $x^*$ as $\deg_{x^*}(w)$, then 
$$H.w = (\deg_{x^*}(w)-\deg_x(w))w.$$
\end{lemma}
\prf
We know that $H.x = -x$ and $H.x^* = x^*$, so
$$w_1\otimes\dots\otimes H.w_i\otimes w_{i+1}\otimes\dots\otimes w_n = (-1)^{\delta_{w_ix}} w_1\otimes\dots\otimes w_n,$$
whence
\begin{eqnarray*}
H.w &=& \sum_{i=1}^n w_1\otimes\dots\otimes h.w_i\otimes w_{i+1}\otimes\dots\otimes w_n\\
&=& \sum_{i=1}^n (-1)^{\delta_{w_ix}} w_1\otimes\dots\otimes w_n \\
&=& (\deg_{x^*}(w)-\deg_x(w))w.
\end{eqnarray*}
\eop
Now let $B$ be the canonical basis of $V^{\otimes n}$ consisting of all words of length $n$ in $x$ and $x^*$, then we can stratify $B$ as
$$B = \bigcup_{m\leq n} B_m$$
with 
$$B_{m} = \{w\in B\mid\deg_x(w) = m,\deg_{x^*}(w) = n-m\}.$$
The number of elements in $B_{m}$ is then equal to
$$b_{m} = {n \choose m},$$
Now Lemma \ref{eigenvector} implies each of these elements is an eigenvector of weight $n-2m$. 
From this, we easily deduce that 
\begin{lemma}\label{wordcount}
In the decomposition of $V^{\otimes n}$ as $\sl_2$-module, the highest weight module with weight $n-2m$ occurs with multiplicity
$${n\choose m} - {n \choose m-1}.$$
That is, if $n$ is odd, this means we have a decomposition
$$V^{\otimes n} = \bigoplus_{m=0}^{\frac{n-1}{2}} V_{n-2m}^{\oplus {n\choose m} - {n \choose m-1}}.$$
If $n$ is even, we have a decomposition
$$V^{\otimes n} = \bigoplus_{m=0}^{\frac{n}{2}} V_{n-2m}^{\oplus {n\choose m} - {n \choose m-1}}.$$
\end{lemma}
\prf
It suffices to observe that the number of elements in $B_m$ counts the total number of eigenvectors of $H$ with weight $n-2m$. The number of elements in $B_{m-1}$ is the total number of eigenvectors with weight $n-2m+2$, each of which is mapped to an eigenvector with weight $n-2m$ under the action of $E$. This means in particular that of all eigenvectors counted in $B_m$, we have exactly $|B_{m-1}|$ eigenvectors which are elements of a highest weight module with weight different from $n-2m$, whence the desired multiplicity.
\eop
In order to obtain a decomposition of the necklace Lie algebra, we have to identify which parts in the decomposition of $V^{\otimes n}$ occur in the decomposition of $C_n$. For this purpose, we first of all need to find a suitable basis for $C_n$. Consider the action of $\ZZ_n$ on $B$ by cyclic permutation. Let $\cO$ be an orbit for this action in $B_n$ and define
$$C(\cO) = \langle w - \overrightarrow{w} \mid w\in\cO\rangle\subset C_n,$$
then the (disjoint) union of all $C(\cO)$ forms a basis for $C_n$, denoted by $\overline{B}$. This basis can again be stratified according to the degree of $\overline{w}$:
$$\overline{B} = \bigcup_{m\leq n} \overline{B}_m$$
with
$$\overline{B}_m = \{\overline{w}\in\overline{B}\mid  \deg_{x^*}(\overline{w}) = m, \deg_{x}(\overline{w}) = n-m\}.$$
Observe that $|C(\cO)| = |\cO|-1$, and the action of $H$ on an element $\overline{w} = w - \overrightarrow{w} \in C(\cO)$ is again with weight $\dim_{x^*}(\overline{w}) - \dim_x(\overline{w})$.
Using this basis we can show
\begin{lemma}\label{necklacecount}
In the decomposition of $C_n$ as $\sl_2$-module, the highest weight module with weight $n-2m$ occurs with multiplicity
$$\sum_{\ell |n | \ell m} \frac{\ell-1}{\ell}\sum_{k |\gcd(\ell,\frac{\ell}{n}m)} \mu(k){\frac{\ell}{k} \choose \frac{\ell m}{kn}}
-
\sum_{\ell |n | \ell m-1} \frac{\ell-1}{\ell}\sum_{k |\gcd(\ell,\frac{\ell}{n}(m-1))} \mu(k){\frac{\ell}{k} \choose \frac{\ell (m-1)}{kn}},
$$
where $\mu$ is the M\"obius function.
\end{lemma}
\prf
Let $\ell | n$, then the Orbit-Stabilizer Theorem tells us each orbit $\cO$ of a point with stabilizer $(\ell)\subset\ZZ_n$ has exactly $\ell$ points, so the number of orbits in $B_m$ of words that have stabilizer $(\ell)\subset\ZZ_n$ is equal to
$$o_\ell = \frac{1}{\ell}|\{w\in B_m\mid \Stab(w) = (\ell)\}|.$$
Next, observe that a word $w$ has stabilizer $(\ell)$ if and only if $w$ can be written up to cyclic permutation as $w= v^{\frac{n}{\ell}}$ with $v$ a word of length $\ell$ with trivial stabilizer and with $\dim_x(v) = \frac{\ell}{n}m$. Observe that in particular this means we only have stabilizer $(\ell)$ if $n|m\ell$.

The orbit of a word of length $\ell$ with trivial stabilizer corresponds to an aperiodic necklace (also known as a Lyndon word) of length $\ell$ with $\frac{\ell}{n}m$ occurrences of $x$. Moreau's necklace counting function tells us this number is given by
$$L_2(\ell,\frac{\ell}{n}m) = \frac{1}{\ell}\sum_{k |\gcd(\ell,\frac{\ell}{n}m)} \mu(k){\frac{\ell}{k} \choose \frac{\ell m}{kn}}.$$
Now we know each orbit contains $\ell$ words, so we have
$$o_\ell = \frac{1}{\ell}\ell L_2(\ell,\frac{\ell}{n}m) = \frac{1}{\ell}\sum_{k |\gcd(\ell,\frac{\ell}{n}m)} \mu(k){\frac{\ell}{k} \choose \frac{\ell m}{kn}}.$$
Each of the orbits $\cO$ counted here gives rise to a basis set $C(\cO)\subset\overline{B}_m$ with $\ell-1$ elements, so for fixed $\ell$, we obtain
$$\frac{\ell-1}{\ell}\sum_{k |\gcd(\ell,\frac{\ell}{n}m)} \mu(k){\frac{\ell}{k} \choose \frac{\ell m}{kn}}$$
distinct basis elements in $\overline{B}_m$. But then the total number of basis elements in $\overline{B}_m$ is simply obtained by summing over all $\ell$ that occur, yielding a total number of elements in $\overline{B}_m$ of
$$\sum_{\ell |n | \ell m} \frac{\ell-1}{\ell}\sum_{k |\gcd(\ell,\frac{\ell}{n}m)} \mu(k){\frac{\ell}{k} \choose \frac{\ell m}{kn}}.$$
Now we conclude by observing again that $|B_{m-1}|$ of these elements come from highest weight modules of weight different from $n-2m$, which concludes the proof.
\eop
Combining Lemma \ref{wordcount} with Lemma \ref{necklacecount} we arrive at the following Clebsch-Gordan decomposition of the necklace Lie algebra.
\begin{theorem}
Let $n\geq 0$ and let $m\leq\frac{n}{2}$, then the highest weight module $V_{n-2m}$ occurs in the decomposition of $\nl_n$ with multiplicity
\begin{align*}
{n\choose m} - {n \choose m-1} &
- \sum_{\ell |n | \ell m} \frac{\ell-1}{\ell}\sum_{k |\gcd(\ell,\frac{\ell}{n}m)} \mu(k){\frac{\ell}{k} \choose \frac{\ell m}{kn}} \\
&+  \sum_{\ell |n | \ell (m-1)} \frac{\ell-1}{\ell}\sum_{k |\gcd(\ell,\frac{\ell}{n}(m-1))} \mu(k){\frac{\ell}{k} \choose \frac{\ell (m-1)}{kn}},
\end{align*}
where $\mu$ is the M\"obius function.
\end{theorem}
In Table \ref{coeff} we list for small total degree the multiplicities of the corresponding highest weight modules.
\begin{table}
\begin{center}
\begin{tabular}{c|c|c|c|c|c|c|c|c|c}
    & 8 & 7 & 6 & 5 & 4 & 3 & 2 & 1 & 0\\
\hline
1 &  0 & 0 & 0 & 0 & 0 & 0 & 0 & 1 & 0\\
2 &  0 & 0 & 0 & 0 & 0 & 0 & 1 & 0 & 0\\
3 &  0 & 0 & 0 & 0 & 0 & 1 & 0 & 0 & 0\\
4 &  0 & 0 & 0 & 0 & 1 & 0 & 0 & 0 & 1\\
5 &  0 & 0 & 0 & 1 & 0 & 0 & 0 & 1 & 0\\
6 &  0 & 0 & 1 & 0 & 0 & 0 & 2 & 0 & 1\\
7 &  0 & 1 & 0 & 0 & 0 & 2 & 0 & 2 & 0\\
8 &  1 & 0 & 0 & 0 & 3 & 0 & 3 & 0 & 3
\end{tabular}
\end{center}
\vspace{.2cm}
\caption{Multiplicities of the highest weight modules in the decomposition of $\nl$. The columns list the weights, whereas the rows list the degrees of the homogeneous parts.}
\label{coeff}
\end{table}

\section{Non-Symplectic Necklace Lie Algebras}\label{nonsymplectic}
Before continuing our study of the symplectic necklace Lie algebra, we would like to make some remarks about generalizations of the results above to non-symplectic necklace Lie algebras, induced by linear double Poisson brackets.

In \cite[Prop. 10]{PVDW}, a one-to-one correspondence was established between finite dimensional associative algebra structures on a vector space $V$ on one hand and linear double Poisson structures on the tensor algebra $T(V)$ on the other hand. Explicitly, let $V$ be a vector space of dimension $n$, generated by elements $x_1, \dots, x_n$, and let $a_{ij}^k$ with $1\leq i,j,k\leq n$ be the structure constants of an associative multiplication defined on $V$, then the bracket
$$\db{x_i,x_j} := \sum_{k=1}^n a_{ij}^k x_k\otimes 1 - a_{ji}^k 1 \otimes x_k$$
is a double Poisson bracket. Moreover, these are the only linear double Poisson brackets on $T(V)$. This leads to the following necklace Lie algebras
\begin{prop}\label{commbrack}
Let $V$ be a finite dimensional associative algebra of dimension $n$, generated by elements $x_i$, $1\leq i \leq n$ and with multiplication
$$x_ix_j = \sum_{k=1}^n a_{ij}^k x_k$$
with $a_{ij}^k\in\CC$ for all $1\leq i,j,k\leq n$. Let $A = T(V) = \CC\langle x_1,\dots,x_n\rangle$ be the tensor algebra of $V$ with double Poisson bracket
$$\db{x_i,x_j} := x_ix_j\otimes 1 - 1 \otimes x_jx_i,$$
then the associated necklace Lie algebra $\nl_A$ has a bracket of degree $-1$ with respect to the natural grading on $T(V)$, for which the degree $1$ component is a Lie subalgebra, equal as a Lie algebra to $(V, [-,-])$. Here the latter bracket is the commutator bracket in the associative algebra $V$:
$$[x_i,x_j] := x_ix_j-x_jx_i.$$
\end{prop}
\prf
This follows immediately from the definition of the double Poisson bracket and the fact that
$$[x_i, x_j] = \sum_{i=1}^k (c_{ij}^k-c_{ji}^k)x_k.$$
\eop

An example of particular importance here would be $V = \mathfrak{gl}_n$ as a finite dimensional associative algebra. The corresponding necklace Lie algebra, which we denote by $\mathfrak{ngl}_n$, would be the necklace Lie algebra of necklaces in $n^2$ beads $e_{ij}$, $1\leq i,j\leq n$, with the bracket induced from the double bracket defined as
$$\db{e_{ij},e_{k\ell}} := \delta_{jk}e_{i\ell}\otimes 1 - \delta_{i\ell}1 \otimes e_{kj}$$
for $1\leq i,j,k,\ell$. The morphism $\varphi$ in Theorem \ref{LieIso} in this situation becomes a Lie algebra morphism from $\mathfrak{ngl}_n$ to $S(\mathfrak{gl}_n)$, inducing an equality as Lie algebras between the degree $1$ component of $\mathfrak{nl}_n$ and the Lie algebra $\mathfrak{gl}_n$ we started with. Moreover, we now have that all homogeneous components of the necklace Lie algebra are Lie algebra modules for the degree $1$ part because of the degree of the bracket. 

This has two important consequences. Firstly, we can again try determine a decomposition of these homogeneous components into simple modules for the degree $1$ part. Secondly, we can construct the necklace equivalents of the classical Lie subalgebras of $\mathfrak{gl}_n$ by restricting the beads of $\mathfrak{ngl}_n$ to the generators of these Lie subalgebras. Suppose $\mathfrak{g}$ is such a Lie subalgebra, and let $\mathfrak{ng}$ be the corresponding Lie subalgebra of $\mathfrak{nl}_n$, then it is easy to see $\mathfrak{ng}$ is a necklace Lie algebra of the tensor algebra $A_\mathfrak{g} := T(\mathfrak{g})$ and hence induces a Poisson structure on all spaces $\iss_d(A_\mathfrak{g})$. One could say that in this way we have defined a \emph{noncommutative thickening} of the Lie-Poisson algebra $S(\mathfrak{g})$.

Observe moreover that the necklace Lie algebras $\mathfrak{ng}$ we have constructed in this way no longer come from double Poisson structures on the tensor algebra $A_\mathfrak{g}$. Indeed, consider for instance the case of $\sl_2$ in $\mathfrak{gl}_2$. To construct $\mathfrak{nsl}_2$, we restrict ourselves to the Lie subalgebra of $\mathfrak{ngl}_2$ on the beads $E = e_{12}$, $F = e_{21}$ and $H = e_{11}-e_{22}$. The bracket $\{E,F\} = H$ is induced by the double Poisson bracket on $T(\mathfrak{gl}_2)$ via
$$\db{e_{12},e_{21}} = e_{11}\otimes 1 - 1 \otimes e_{22},$$
which is not in $A_{\mathfrak{nsl}_2}\otimes A_{\mathfrak{nsl}_2}$ when considering $A_{\mathfrak{nsl}_2}$ as embedded in $T(\mathfrak{gl}_2)$. As a matter of fact, it follows from \cite[Prop. 10]{PVDW} and Proposition \ref{commbrack} that the necklace Lie bracket constructed this way cannot be induced from a double Poisson bracket on $A_{\mathfrak{nsl}_2}$ as the Lie bracket of $\sl_2$ is not induced as the commutator bracket of an associative multiplication on $\sl_2$.

\section{The Center of the Necklace Lie Algebra}\label{center}
From Corollary \ref{notsimplealgebra} we already knew that the necklace Lie algebra was not a simple Lie algebra. In fact, more is true, and in this section we will show that the necklace Lie algebra in fact has a nontrivial, infinite dimensional center. We will denote the center of $\nl$ by $\mathfrak{c}$. We have
\begin{theorem}
For all $n\in\NN$, the element $c_n = [x,x^*]^n \in\nl$ is a central element. These elements are distinct and nontrivial for $n>1$, and hence $\mathfrak{c}$ is infinite dimensional.
\end{theorem}
\prf
A straightforward computation on the level of the double Poisson algebra $\CC\langle x,x^*\rangle$ yields
\begin{eqnarray*}
\db{[x,x^*],x} &=& x\otimes 1 - 1 \otimes x \\
\db{[x,x^*],y} &=& y\otimes 1 - 1 \otimes y
\end{eqnarray*}
Which yields because of the twisted antisymmetry and the derivation property of a double Poisson bracket that
\begin{eqnarray*}
\db{[x,x^*]^n,x} &=& \sum_{i=0}^{n-1} (x[x,y]^{n-1-i}\otimes [x,y]^{i} - [x,y]^{n-1-i}\otimes [x,y]^{i}x) \\
\db{[x,x^*]^n,x^*} &=& \sum_{i=0}^{n-1} (x^*[x,x^*]^{n-1-i}\otimes [x,x^*]^{i} - [x,x^*]^{n-1-i}\otimes [x,y]^{i}y).
\end{eqnarray*}
From which we deduce that for the Poisson-Loday bracket we have
\begin{eqnarray*}
\{[x,x^*]^n,x\}_L &=& n(x[x,x^*]^{n-1} - [x,x^*]^{n-1}x) \\
\{[x,x^*]^n,x^*\}_L &=& n(x^*[x,x^*]^{n-1} - [x,x^*]^{n-1}x^*),
\end{eqnarray*}
and so
\begin{eqnarray*}
\{[x,x^*]^n,x^a(x^*)^b\}_L &=& n(x^a(x^*)^b[x,x^*]^{n-1}-[x,x^*]^{n-1}x^a(x^*)^b).
\end{eqnarray*}
But then for a general word $w= x^{a_1}(x^*)^{b_1}\dots x^{a_k}(x^*)^{b_k}$ we obtain
\begin{eqnarray*}
\{[x,x^*]^n,w\}_L &=& \sum_{i=1}^k
\left (\prod_{j=1}^{i-1} x^{a_j}(x^*)^{b_j}\right )
\{x,x^{a_i}(x^*)^{b_i}\}_L
\left (\prod_{\ell=i+1}^{k} x^{a_\ell}(x^*)^{b_\ell}\right )\\
&=& n\sum_{i=1}^k
\left (\prod_{j=1}^{i-1} x^{a_j}(x^*)^{b_j}\right )
(x^a(x^*)^b[x,x^*]^{n-1}\\
& &-[x,x^*]^{n-1}x^a(x^*)^b)
\left (\prod_{\ell=i+1}^{k} x^{a_\ell}(x^*)^{b_\ell}\right )\\
&=& -n([x,y]^{n-1}w - w[x,y]^{n-1}),
\end{eqnarray*}
which is easily seen to be zero modulo commutators and hence 
$$\{[x,x^*]^n,w\} = 0 \in \nl.$$

To establish that these elements are all distinct, it suffices to evaluate $c_n$, for $n>1$, in the matrix couple $(X,X^*)$ given by
$$X = \begin{pmatrix} 0 & \lambda & 0 \\ 0 & 0 & -\lambda \\ 0& 0 & 0 \end{pmatrix}
\quad\mathrm{and}\quad
X^* = \begin{pmatrix}0 & 0 & 0 \\ 1 & 0 &0 \\ 0 & 1 & 0\end{pmatrix},$$
with $\lambda\in\CC^*$. Indeed, we have that 
$$c_n(X,X^*) = tr(\begin{pmatrix} \lambda & 0 & 0 \\
0 & -2\lambda & 0 \\ 0 & 0 & \lambda \end{pmatrix}^n)
= 2\lambda^n + (-2)^n\lambda^n \neq 0.
$$
\eop

\begin{remark}
Observe that the description of the center obtained here is actually a special case of the stronger result \cite[Thm. 8.6.1, Cor. 8.6.2]{CBEG}. However, we feel that the proof presented here, through direct computation using the double Poisson structure, provides an interesting alternative approach for the particular case considered in this note.
\end{remark}

An immediate consequence for the induced Poisson geometry on the quotient spaces is
\begin{coro}
For any $n\geq 2$, the zeroth Poisson cohomology 
$$PH^0(\CC[\iss_n(\CC\langle x,x^*\rangle)])$$
is infinite dimensional.
\end{coro}

\section{The Poisson Geometry of $C_2$}\label{leaves}
As an extended example of the Poisson geometry induced by the necklace Lie algebra, we will have a look at the Poisson structure induced by $\nl$ on small dimensional representation spaces. First of all observe that the space of one-dimensional representations, 
$$\rep_1(\CC\langle x,x^*\rangle) = \iss_1(\CC\langle x,x^*\rangle),$$
simply becomes the symplectic plane when equipped with the induced Poisson structure. 

The first degenerate case we encounter already is the variety 
$$\iss_2(\CC\langle x,x^*\rangle).$$
It is known that this is a five dimensional affine space, and so cannot be a symplectic manifold. The coordinate ring of this space is generated as a polynomial algebra by $tr(x)$, $tr(x^*)$, $tr(x^2)$, $tr((x^*)^2)$ and $tr(xx^*)$. 

The Poisson structure induced on the coordinate ring 
$$C = \CC[\iss_2(\CC\langle x,x^*\rangle)]$$
then exactly corresponds to the Poisson structure induced by $\nl_{\leq 2} = \sl_2\rtimes\h$ on the quotient of its symmetric algebra (the center of which was studied in \cite{AOV}) by the ideal $(tr(1) - 2)$. We summarized this structure in the multiplication table listed in Table \ref{PoissonStructureTable}.
\begin{table}
\begin{center}
\begin{tabular}{c|c|c|c|c|c}
$\{-,-\}$         & $tr(x)$ & $tr(x^*)$ & $tr(x^2)$ & $tr((x^*)^2)$ & $tr(xx^*)$ \\
\hline
$tr(x)$          &   $0$    &    $2$      &     $0$     &  $2tr(x^*)$    & $tr(x)$      \\
$tr(x^*)$       &   $-2$   &    $0$      & $-2tr(x)$ &   $0$             & $-tr(x^*)$ \\
$tr(x^2)$      &   $0$    & $2tr(x)$   & $0$         & $4tr(xx^*)$   & $2tr(x^2)$ \\
$tr((x^*)^2)$ & $-2tr((x^*)^2$ & $0$ & $-4tr(xx^*)$ & $0$        & $-2tr((x^*)^2$ \\
$tr(xx^*)$     & $-tr(x)$ & $tr(x^*)$  & $-2tr(x^2)$   & $2tr((x^*)^2)$ & 0
\end{tabular}
\end{center}
\vspace{.2cm}
\caption{The Poisson structure on the generators of $C$.}
\label{PoissonStructureTable}
\end{table}
Denote the canonical generators of $\sl_2$ by $E$, $F$ and $H$ and denote the canonical generators of $\h$ by $X$, $Y$ and $Z$ then Table \ref{PoissonStructureTable} yields the identification $H = tr(xy)$, $F= - \frac{tr(x^2)}{2}$, $E = \frac{tr((x^*)^2)}{2}$, $X = tr(x^*)$, $Y = -tr(x)$ and $Z = 2$.
Now consider the following automorphism of $\frac{S(\sl_2\rtimes\h)}{(Z-2)}$:
\begin{eqnarray*}
H' &=& H + \frac{XY}{2}\\
E' &=& E - \frac{X^2}{4}\\
F' &=& F + \frac{Y^2}{4}\\
X' &=& X\\
Y' &=& Y\\
\end{eqnarray*}
Then we obtain that 
$$\frac{S(\sl_2\rtimes\h)}{(Z-2)} \cong \frac{S(\sl_2)\otimes S(\h)}{(Z-2)}$$
as Poisson algebras. Indeed, observe for instance that
\begin{eqnarray*}
\{H',E'\} &=& \{H,E\} - \{H, \frac{X^2}{4}\}+ \{\frac{XY}{2},E\} - \{\frac{XY}{2},\frac{X^2}{4}\}\\
&=& 2E - \frac{X^2}{2} -\frac{X^2}{2} + \frac{X^2Z}{4}\\
&=& 2E'.
\end{eqnarray*}
Similarly one obtains $\{H',F'\} = -2F'$, $\{E',F'\} = H'$, identifying the subalgebra generated by $H'$, $E'$ and $F'$ with $\sl_2$. Analogously, one observes that $H'$, $E'$ and $F'$ act trivially on $X'$ and $Y'$, finishing the claim.

%\begin{remark}
%It is interesting to remark that we have actually shown that we have an isomorphism
%$$S(\sl_2\rtimes\h) \cong S(\sl_2)\otimes S(\h),$$
%whereas the Lie algebras $\sl_2\rtimes\h$ and $\sl_2\times \h$ are definitely non-isomorphic as the former is a perfect Lie algebra and the latter is not.
%\end{remark}
\begin{remark}
When dealing with $n\times n$ matrix invariants, a canonical way of obtaining nicer relations is by looking at traceless matrices. That is, instead of looking at the generators $tr(w)$ where $w$ is a word in the generic matrices, one looks at the generic matrices $tr(x_i)$ in combination with $tr(w')$ where $w'$ is a word in the \emph{traceless generic matrices} $x_i - \frac{1}{n}tr(x_i)$. It is a general principle that the relations with respect to these generators behave more nicely.

Now observe that the automorphism above actually corresponds to changing from traces of generic matrices to traces of words in traceless matrices. For instance,
$$H' = tr(xx^*)-\frac{tr(x)tr(y)}{2} = tr((x-\frac{1}{2}tr(x))(x^*-\frac{1}{2}tr(x^*))).$$
\end{remark}

But now 
$$\frac{S(\sl_2)\otimes S(\h)}{(Z-2)}\cong S(\sl_2)\otimes \frac{S(\h)}{(Z-2)},$$
which means in particular that the center of $\frac{S(\sl_2)\otimes S(\h)}{(Z-2)}$ is equal to the center of $S(\sl_2)$. This center is known to be the polynomial algebra generated by the Casimir
$$c_{\sl_2} = {H'}^2 + 4E'F'.$$
This allows us to conclude the Poisson center of $C$ is generated as a polynomial algebra by the element $c_{\sl_2}$ under the identification listed above. We have
\begin{prop}
The canonical morphism $S(\nl)\twoheadrightarrow C$  maps for $n\geq 0$
\begin{eqnarray*}
c_{2n} &\mapsto & -2^{1-n}c_{\sl_2}^{n} \\
c_{2n+1} & \mapsto & 0
\end{eqnarray*}
\end{prop}
\prf
This is an immediate consequence of the \emph{Cayley-Hamilton} relations satisfied by matrices. Recall that for any $2\times 2$ matrix $A$ we have
$$A^2 - tr(A)A + det(A) = 0,$$
and
$$det(A) = \frac{1}{2}(tr(A)^2 - tr(A^2)).$$
This means in particular that
\begin{eqnarray*}
tr([x,x^*]^{2n}) &=& tr( (tr([x,x^*])[x,x^*] - \frac{1}{2}(tr([x,x^*])^2 - tr([x,x^*]^2)))^n) \\
&=& tr((\frac{1}{2}tr([x,x^*]^2))^n) \\
&=& 2^{1-n}tr([x,x^*]^2)^n,
\end{eqnarray*}
and
\begin{eqnarray*}
tr([x,x^*]^{2n+1}) &=& tr( (tr([x,x^*])[x,x^*] - \frac{1}{2}(tr([x,x^*])^2 - \\
&  &tr([x,x^*]^2)))^n[x,x^*]) \\
&=& tr((\frac{1}{2}tr([x,x^*]^2))^n[x,x^*]) \\
&=& 2^{1-n}tr([x,x^*]^2)^ntr([x,x^*]) \\
&=& 0.
\end{eqnarray*}
Now we compute
\begin{eqnarray*}
tr([x,x^*]^2) &=& tr(x^2(x^*)^2) - tr((xx^*)^2) \\
&=& tr( (tr(x)x - det(x))(tr(x^*)x^*-  det(x^*))) \\
&  &- tr(tr(xx^*)xx^*-det(x)det(x^*))\\
&=& tr(x)tr(x^*)tr(xx^*) - det(x)tr(x^*)^2 - det(x^*)tr(x)^2 \\
&  & + 2det(x)det(x^*)- tr(xx^*)^2 + 2det(x)det(x^*) \\
&=& tr(x)tr(x^*)tr(xx^*) - tr(xx^*)^2  + tr(x^2)tr((x^*)^2)\\
&&  - \frac{1}{2}(tr(x^2)tr(x^*)^2 + tr(x)^2tr((x^*)^2),
\end{eqnarray*}
which finishes the proof under the identification $H = tr(xy)$, $F= - \frac{tr(x^2)}{2}$, $E = \frac{tr((x^*)^2)}{2}$, $X = tr(x^*)$ and $Y = -tr(x)$.
\eop

Another consequence of the isomorphism above is the following.
\begin{theorem}
The symplectic leaves of $\iss_2(\CC\langle x,x^*\rangle)$ correspond to the coadjoint orbits of $\sl_2$ and hence correspond to the level sets of the Casimir element $c_{\sl_2}$
$$\mathcal{S}_\lambda := \{\rho\in\iss_2(\CC\langle x,x^*\rangle)\mid c_{\sl_2}(\rho) = \lambda\},$$
for $\lambda\in\CC^*$ together with the two sets
$$\mathcal{S}_0' := \{\rho\in\iss_2(\CC\langle x,x^*\rangle)\mid c_{\sl_2}(\rho) = 0,
(E',F',H') \neq 0\},$$
and
$$\mathcal{S}_0'' := \{\rho\in\iss_2(\CC\langle x,x^*\rangle)\mid (E',F',H') = 0\},$$
\end{theorem}
\prf
This follows immediately from the above coordinate change, the fact that the Poisson structure on $\frac{S(\h)}{(Z-2)}$ is the canonical symplectic structure and the fact that the Poisson structure on $S(\sl_2)$ is a Lie-Poisson structure, the symplectic leaves of which correspond to the coadjoint orbits. These are known to be equal to the description in the theorem.
\eop

We would like to point out a relationship between the symplectic leaves and the \emph{Luna stratification} of $\iss_2(\CC\langle x,x^*\rangle)$. Recall for a general associative algebra $A$, elements in $\iss_n(A)$ correspond to isomorphism classes of semi simple representations. That is, an element $\rho\in \iss_n(A)$ corresponds to an isomorphism class of
$$\sigma_1^{\oplus e_1}\oplus\dots\oplus\sigma_k^{\oplus e_k},$$
with all $\sigma_i$ distinct simple representations of dimension $s_i$ and
$s_1e_1+\dots+s_ne_n = 0.$ The unordered $k$-tuple $\tau = [(s_1,e_1);\dots;(s_k,e_k)]$ is called the \emph{representation type} of $\rho$, and for fixed $\tau$ the set $S_\tau$ of all elements in $\iss_n(A)$ of representation type $\tau$ is called the \emph{Luna stratum of representation type $\tau$}.

For $\iss_2(\CC\langle x,x^*\rangle)$ we have three Luna strata corresponding to representation type $\tau_1 = [(2,1)]$, $\tau_2 = [(1,1);(1,1)]$ and $\tau_3 = [(1,2)]$. Using the identification from the beginning of this section, consider
\begin{eqnarray*}
\mathcal{H} &:=& \{ (X,Y,E,F,H)\in\iss_2(\CC\langle x,x^*\rangle) \mid c_{sl_2}(X,Y,E,F,H) = 0\} \\
\mathcal{S} &:=& \{ (X,Y,E,F,H)\in\iss_2(\CC\langle x,x^*\rangle) \mid E' = F' = H' = 0\},
\end{eqnarray*}
then (for the actual computations we refer to \cite{LievenBook})
\begin{eqnarray*}
S_{\tau_1} &=& \iss_2(\CC\langle x,x^*\rangle)\backslash \mathcal{H}\\
S_{\tau_2} &=& \mathcal{H}\backslash\mathcal{S} \\
S_{\tau_3} &=& \mathcal{S}.
\end{eqnarray*}
So observe that the stratification by symplectic leaves refines the simple stratum, whereas the other strata each are exactly one coadjoint orbit.

\begin{remark}
Observe that because in this setting the Poisson structure on $\iss_n(A)$ is induced by a Poisson structure on $\rep_n(A)$, the stratification by symplectic leaves will always be finer than the Luna stratification (see \cite[Prop. 10.5.2]{Marsden} or \cite[Lemma 4.3]{Martino}). This means that checking whether or not a given representation is simple becomes equivalent to checking whether or not the representation belongs to a certain family of symplectic leaves. As the Poisson structure induced on $\iss_n(A)$ actually comes from a linear double Poisson structure, one could hope that these symplectic leaves still have a nice description, thus providing an easy method for checking whether a given representation is simple or not. Moreover, the symplectic flows give a canonical way of constructing new, non-isomorphic simple representations from a given representation. 
\end{remark}

\section{Poisson Orders and the Trace Algebra}\label{order}
In this final section, $A$ is once again any affine associative algebra.
We begin with some definitions, motivated by the definition of a formal necklace algebra and a formal trace algebra in \cite{LievenBook}.
\begin{defn}
Let $\frac{A}{[A,A]}$ be the vector space of necklaces of an algebra $A$, then the \emph{necklace algebra of $A$} is the symmetric algebra
$$\cN := S(\frac{A}{[A,A]})$$
generated by all necklaces. The \emph{trace algebra of $A$} is the algebra
$$\cT := \cN\otimes_\CC A.$$
\end{defn}
Now recall the notion of a Poisson order, introduced in \cite{PoissonOrder}.
\begin{defn}
Let $A$ be an algebra with a central subalgebra $Z$, and let $H : Z\rightarrow\Der_\CC(A)$ be a linear map satisfying the Leibniz identity, then $(A,Z,H)$ is called a \emph{Poisson order} if
\begin{enumerate}
\item $A$ is finitely generated as a module over $Z$;
\item $H\mid_Z$ defines a Poisson structure on $Z$.
\end{enumerate}
\end{defn}
We will need a slightly relaxed version of this definition.
\begin{defn}
Let $A$, $Z$ and $H$ be as in the previous definition, then $(A,Z,H)$ is called a \emph{formal Poisson order} if $H\mid_Z$ defines a Poisson structure on $Z$.
\end{defn}

We have
\begin{prop}
Let $A$ be a double Poisson algebra with bracket $\db{-,-}$ and corresponding necklace Lie algebra $\nl$. Consider the Lie-Poisson bracket $\{-,-\}$ induced by $\db{-,-}$ on $\cN$ and define
$$\tilde{H}: \cN\rightarrow \cN : w \mapsto \{w,-\},$$
Extend this map using the Leibniz rule to a map $H$ on $\cT$,
then $(\cT,\cN,H)$ is a formal Poisson order.
\end{prop}
\prf
This follows immediately from the well-definedness of the map $\tilde{H}$ and the definition of a formal Poisson order.
\eop
Which implies
\begin{coro}
Let $\cT_n^d$ be the trace ring of $d$-tuples of generic $n\times n$ matrices, then any double Poisson bracket on $\CC\langle x_1,\dots,x_d\rangle$ turns $\cT_n^d$ into a Poisson order over its center.
\end{coro}
\prf
This follows from the fact that the map $\tilde{H}$ is equivariant under taking traces and hence respects the formal Cayley-Hamilton relations which form the kernel of the canonical map $\cT\twoheadrightarrow \cT_n^d$.
\eop

\begin{remark}
It was shown in \cite{PoissonOrder} that the stratification by \emph{symplectic cores} has a nice link with the representation theory of a Poisson order. More precisely, given a Poisson order $(A,Z,H)$ and a point $x\in MaxSpec(Z)$, one can look at the finite dimensional algebra $A_x := \frac{A}{\mathfrak{m}_xA}$, with $\mathfrak{m}_x$ the defining maximal ideal of $x$. Now \cite[Theorem 4.2]{PoissonOrder} shows that if $x$ and $y$ belong to the same core, $A_x$ and $A_y$ are isomorphic. An interesting question one can ask now is whether or not it is always possible to find a double Poisson bracket such that this result is optimal for the trace ring of the double Poisson algebra. That is, can one find a double Poisson bracket such that $A_x$ and $A_y$ are isomorphic \emph{if and only if} $x$ and $y$ belong to the same core?
\end{remark}

\end{document}